\documentclass[10pt,twoside]{article}

\usepackage{times, fancyhdr}
\usepackage{cite}
\usepackage[ ]{graphicx}
\usepackage{epsfig}
\usepackage{amsmath}
\usepackage{amssymb}
\usepackage{amsfonts}

\usepackage{amsthm}
\usepackage{indentfirst}
\usepackage{sectsty}

\usepackage{float}

\numberwithin{equation}{section}

\begin{document}

\title{\large Existence of double  Walsh series universal in weighted
$L_\mu^1[0,1]^2$ spaces} 
\date{}

\author{\normalsize S.A.Episkoposian
\\ {\footnotesize Received May 25, 2007; Revised June 19, 2007} }
\maketitle

\thispagestyle{empty}

\begin{abstract}
\noindent In this paper we consider a question on existence of
double  Walsh series universal in weighted $L_\mu^1[0,1]^2$
spaces. We construct a weighted function $\mu(x,y)$ and a series
by double Walsh system of the form
$$\sum_{n,k=1}^\infty c_{n,k}W_n(x)W_k(y)\ \  \mbox{with} \ \ \sum_{n,k=1}^\infty \left |
c_{n,k} \right|^q <\infty\ \mbox{for all}\ q>2,$$ which is universal in
$L_\mu^1[0,1]^2$ concerning subseries with respect to convergence,
in the sense of both spherical and rectangular partial sums.
\smallskip
\smallskip
\smallskip



\end{abstract}

\section{Introduction}

Let $\mu(x)$, $0<\mu(x) \le1, x\in[0,1]$
 be a measurable on $[0,1]$ function and let $L_\mu^1[0,1]$ be a space of real measurable functions
 $f(x)$, $x\in [0,1]$ with

\[
  \int_0^1 |f(x)| \mu(x) dx<\infty.
  \]
{\bf Definition 1.} A functional series
\begin{equation}
 \sum_{k=1}^\infty f_k(x),\ \ f_k(x) \in L_\mu^1[0,1]
\end{equation}
is called universal in weighted space $L_\mu^1[0,1]$ with respect
to  rearrangements, if for any function $f(x) \in L_\mu^1[0,1]$
the terms of (1.1) can be rearranged so that the obtained series
$\displaystyle {\sum_{k=1}^\infty f_{\sigma(k)}(x)}$ converges to
$f(x)$ in the metric $L_\mu^1[0,1]$,  i.e.

\[
\lim_{n\to \infty} \int_0^{1} \left| \sum _{k=1}^n f_{ \sigma
(k)}(x)-f(x) \right| \cdot \mu(x) dx=0.
\]
{\bf Definition 2.} The series (1.1) is called universal in
$L_\mu^1[0,1]$ concerning subseries, if for any function $f(x) \in
L_\mu^1[0,1]$ there exists a subseries $\displaystyle
{\sum_{k=1}^\infty f_{n_k}(x)}$ of (1.1), which converges to
$f(x)$ in the metric $L_\mu^1[0,1]$. \vskip 4mm

The above mentioned definitions are  given not in the most general
form but only in the form which will be applied in the present
paper. Note, that for one-dimensional case there are many papers
are devoted to the question on existence of various types of
universal series in the sense of convergence almost everywhere and
on a measure ( see [1], [3]- [9], [11], [12]).

The first usual universal in the sense of convergence almost
everywhere trigonometric series were constructed by D.E.Menshov
[8] and V.Ya.Kozlov [7]. The series of the form
 $${1\over2}+\sum_{k=1}^\infty a_k \cos{kx}+b_k \sin{kx} \eqno(A)$$
was constructed just by them such that for any measurable on
$[0,2\pi]$  function $f(x)$ there exists the growing sequence of
natural numbers $n_k$  such that the series (A) having the
sequence of partial sums with numbers $n_k$   converges to $f(x)$
almost everywhere on $[0,2\pi]$.

Note that in this result, when $f(x)\in{L^1_{[0,2\pi]}} $, it is
impossible to replace convergence almost everywhere by convergence
in the metric ${L^1_{[0,2\pi]}}$.

 This result was distributed by A.A.Talalian on arbitrary orthonormal
 complete systems (see [11]). He also established (see [12]), that if
 $\{\phi_n(x)\}_{n=1}^\infty $  - the normalized basis of space ${L^p_{[0,1]}},p>1 $,
  then there exists a series of the form
  $$\sum_{k=1}^\infty{a_k\phi_k(x)},\ \ a_k \to 0.\eqno(B)$$
which has property: for any measurable function $f(x)$ the members
of series (B) can be rearranged so that the again received series
converge on a measure on [0,1] to $f(x)$.

In [2] these results are transferred to two-dimensional case.

W. Orlicz [9] observed the fact that there exist functional series
that are universal with respect to rearrangements in the sense of
a.e. convergence in the class of a.e. finite measurable functions.
It is also useful to note that even Rieman proved that every
convergent numerical series which is not absolutely convergent is
universal with respect to rearrangements in the class of all real
numbers. In [4] and [3] the following results are proved.

{\bf Theorem 1.} {\it There exists a series of the form
\begin{equation}
 \sum_{k=1}^\infty{c_kW_k(x) \ \  with \ \ \sum_{k=1}^\infty \left | {c_k} \right|^q <\infty}\ \ \ for\ all\ q>2
\end{equation}
such that for any number $\varepsilon>0$  a weighted function
$\mu(x)$ with
\begin{equation}
 0<\mu(x) \le1, \left | \{ x\in[0,1]: \mu(x)\not =1
\} \right | <\varepsilon
\end{equation}
can be constructed, so that the series (1.2) is universal in
$L_\mu^1[0,1]$ with respect to rearrangements.}

 {\bf Theorem 2.} {\it There exists a series of the form (1.2)
such that for any number $\varepsilon>0$  a weighted function
$\mu(x)$ with (1.3) can be constructed, so that the series (1.2)
is universal in $L_\mu^1[0,1]$ concerning subseries.}

In this paper we prove that Theorems 1 and 2 can be transferred
from one-dimensional case to two-dimensional one.

Moreover, the following statements are true.

{\bf Theorem 3.} {\it There exists a double series of the form
\begin{equation}
\sum_{n,k=1}^\infty{c_{n,k}W_n(x)W_k(y) \ \  with \ \
\sum_{n,k=1}^\infty \left | {c_{n,k}} \right|^q <\infty} \ \ for \
all\ q>2
\end{equation}
with the following property: \noindent for any number $\varepsilon>0$
a weighted function $\mu(x,y)$ satisfying
\begin{equation}
0<\mu(x,y) \le 1, \left | \{ (x,y)\in T=[0,1]^2: \mu(x,y)\not =1
\} \right | <\varepsilon
\end{equation}
  can be constructed so that the series (1.4) is universal
  in $L_\mu^1(T)$ concerning subseries
with respect to convergence in the sense of both spherical and
rectangular partial sums.}

{\bf Theorem 4.} {\it There exists a double series of the form
(1.4) with the following property: \noindent for any number
$\varepsilon>0$ a weighted function $\mu(x,y)$ with (1.5) can be
constructed, so that the series (1.4) is universal in $L_\mu^1(T)$
concerning rearrangements  with respect to convergence in the
sense of both spherical and rectangular partial sums.}

{\bf Remark.} We can also prove that Theorems 3 and 4 remain true
for double trigonometric system. \vskip 4mm

The author thanks Professor M.G.Grigorian for his attention to
this paper.


\section{Notations }

\vskip 4mm

 First we will give a definition of one-dimensional Walsh-Paley system
(see [10]).
\begin{equation}
W_0(x)=1,\ \  W_n(x)= \prod_{s=1}^n r_{m_s}(x),\ \ n= \sum_{s=1}^k
2^{m_s},\ \ m_1>m_2>...>m_s
\end{equation}
where $\{ r_k(x)\}_{k=0}^\infty$ is the system of Rademacher given
as follows:
\[
r_0(x)= 1\ \ \mbox{for}\ \ x \in \bigg[0,{1\over2}\bigg), \ \
r_0(x)=-1 \ \mbox{for}\ \ x \in \bigg({1\over2},1\bigg];
\]
\[
r_0(x+1)=r_0(x)\ \ \mbox{and}\ \ r_k(x)=r_0(2^kx)\ \mbox{for}\ k=1,2,...
\]

The rectangular and spherical partial sums of the double series
\[
\sum_{k,\nu=1}^\infty c_{k,\nu} W_k(x)W_\nu(y)
\]
will be denoted by
\[
 S_{n,m}(x,y)= \sum _{k=1}^n\sum _{\nu=1}^m c_{k,\nu} W_k(x)W_\nu(y)
\]
and
\[
 S_R(x,y)= \sum _{\nu^2+k^2 \leq R^2} c_{k,\nu} W_k(x)W_\nu(y).
\]

If $g(x,y)$ is a continuous function on $T=[0,1]^2$, then we set
 \[
 ||g(x,y)||_C=\max_{(x,y)\in T} |g(x,y)|.
 \]
For an interval $\Delta$ in the form
$$\Delta_m^{(i)}= \left[ {{i-1}\over
{2^m}},{i\over {2^m}} \right],\ \  1\leq i \leq 2^m.$$
we let  $\chi_{\Delta}(x)$ denote
the characteristic function of $\Delta$.

\section{Some basic lemmas}

\vskip 4mm In [3] we proved the following lemma.

\vskip 4mm

{\bf Lemma 1.} {\it For any given numbers $0<\varepsilon<1$,
$N_0>2$ and a step function
$$f(x)= \sum_{s=1}^q \gamma_s \cdot \chi_{\Delta_s} (x),$$
where  $\Delta_s$ is an interval  of the form $\displaystyle
\Delta_m^{(i)}= \left[ {{i-1}\over {2^m}},{i\over {2^m}} \right]
$, $ 1\leq i \leq 2^m $, there exist a measurable set $E \subset
[0,1]$ and a polynomial $P(x)$ of the form
 $$P(x)= \sum_{k=N_0}^N c_kW_k(x) $$
which satisfy the conditions:
$$P(x)=f(x)\ \  on \ \ E, \leqno(1)$$
$$|E|> (1- \varepsilon), \leqno(2)$$
$$\sum_{k=N_0}^N |c_k|^{2+ \varepsilon}< \varepsilon, \leqno(3)$$
$$\max_{N_0 \leq m<N} \left[ \int_e \left | \sum_{k=N_0}^m c_k
W_k(x) \right | dx \right] <\varepsilon+\int_e |f_(x)|dx,\leqno(4)$$
for every measurable subset $e$  of $E$.}
\par\par\bigskip

The following is true .

{\bf Lemma 2.} {\it For any numbers $\gamma \not=0$, $0<\delta<1$,
$N>1$ and for any square $\Delta=\Delta_1 \times \Delta_2 \subset
T=[0,1]^2$ there exists a measurable set $E \subset T$ and a
polynomial $P(x,y)$ of the form
$$P(x,y)= \sum_{k,s=N}^M c_{k,s}W_k(x) \cdot W_s(y) ,$$
with the following properties:
$$|E|> 1- \delta, \leqno(1)$$
$$\sum_{k,s=N}^M |c_{k,s}|^{2+\delta}< \delta, \leqno(2)$$
$$P(x,y)=\gamma \cdot \chi_\Delta (x,y) \ \  for \ \ (x,y) \in E, \leqno(3)$$
$$\max_{N \leq \overline {n}, \overline{m} \leq M} \left[\int \int_e \left | \sum_{k,s=N}^{\overline {n}, \overline{m}}  c_{k,s}W_k(x) \cdot W_s(y) \right | dx dy \right]\leqno(4)$$
$$+\max_{\sqrt {2}N \leq R \leq \sqrt {2}M} \left[ \int\int_e \left | \sum_{2N^2 \leq k^2+s^2 \leq R^2}  c_{k,s}W_k(x) \cdot W_s(y) \right | dx dy \right] \leq 16 \cdot |\gamma|\cdot |\Delta|,$$
for every measurable subset $e$  of $E$.}

{\bf Proof of Lemma 2.} We apply Lemma 1, setting
$$f(x)= \gamma \cdot \chi_{\Delta_1} (x),\ \ N_0=N,\ \ \varepsilon= {\delta \over 2}.$$
Then we can define a measurable set $E_1 \subset [0,1]$ and a
polynomial $P_1(x)$ of the form
 $$P_1(x)= \sum_{k=N}^{N_1} a_kW_k(x) $$
which satisfy the conditions:
$$P_1(x)=\gamma \cdot \chi_{\Delta_1} (x) \ \  for \ \ x \in E_1, \leqno(1^0)$$
$$|E_1|> 1- {\delta \over 2}, \leqno(2^0)$$
$$\sum_{k=N}^{N_1} |a_k|^{2+\delta}< \delta, \leqno(3^0)$$
$$\max_{N \leq \overline {n} \leq N_1} \left[ \int_{e_1} \left | \sum_{k=N}^{\overline {n}}  a_kW_k(x)  \right | dx \right] \leq 2 \cdot |\gamma| \cdot |\Delta_1|, \leqno(4^0)$$
for every measurable subset $e_1$  of $E_1$.

Set
\begin{equation}
M_0=2 \cdot \left( N_1^2+1 \right)
\end{equation}
and apply Lemma 1 again, setting
$$f(y)=  \chi_{\Delta_2} (y),\ \ N_0=M_0,\ \ \varepsilon= {\delta \over 2}.$$
Then we can define a measurable set $E_2 \subset [0,1]$ and a
polynomial $P_2(y)$ of the form
 $$P_2(y)= \sum_{s=M_0}^M b_sW_s(y), $$
which satisfy the conditions:
$$P_2(y)= \chi_{\Delta_2} (y) \ \  for \ \ y \in E_2, \leqno(1^{00})$$
$$|E_2|> 1- {\delta \over 2}, \leqno(2^{00})$$
$$\sum_{s=M_0}^M |b_s|^{2+\delta}< \delta, \leqno(3^{00})$$
$$\max_{M_0 \leq \overline {m} \leq M} \left[ \int_{e_2} \left | \sum_{s=M_0}^{\overline {m}}  b_sW_s(y)  \right | dy \right] \leq 2 \cdot  |\Delta_2|, \leqno(4^{00})$$
for every measurable subset $e_2$  of $E_2$.

Set
\begin{equation}
E=E_1 \times E_2,
\end{equation}
\begin{equation}
P(x,y)= P_1(x) \cdot P_2(x)=\sum_{k,s=N}^M c_{k,s}W_k(x) \cdot
W_s(y) ,
\end{equation}
where
\begin{equation}
c_{k,s}=a_k \cdot b_s,\ \ if \ \ N \leq k \leq N_1, \ \ M_0 \leq s
\leq M
\end{equation}
and
\[
  c_{k,s}=0,\ \ for \ \ other \ \ k\ \ and\ \  s.
\]
By $(1^0) - (3^0),\ \ (1^{00}) - (3^{00})$ and (3.2) - (3.4) we
obtain
$$|E|> 1- \delta, $$
$$\sum_{k,s=N}^M |c_{k,s}|^{2+\delta}= \sum_{k=N}^{N_1}
|a_k|^{2+\delta} \cdot \sum_{s=M_0}^M |b_s|^{2+\delta}< \delta, $$
$$P(x,y)=\gamma \cdot \chi_\Delta (x,y) \ \  for \ \ (x,y) \in E.$$
Thus, the statements 1) - 3) of Lemma 2 are satisfied. Now we will
check the fulfillment of statement 4) of Lemma 2.

Let $N^2+M_0^2<R^2<N_1^2+M^2$, then for some $m_0>M_0$ we have
$m_0<R^2<m_0+1$ and from (3.1) it follows, that
$R^2-N_1^2>(m_0-1)^2$.

Consequently taking relations $(4^0), (4^{00})$ and (3.2) - (3.4)
for any measurable set $e \subset E$ $\left( e=e_1  \times e_2,\ \
e_1 \subset E_1,\ \ e_2 \subset E_2 \right)$ we obtain
$$ \int\int_e \left | \sum_{N^2+M^2 \leq
k^2+s^2 \leq R^2}  c_{k,s}W_k(x) \cdot W_s(y) \right | dx dy $$
$$\leq \int\int_e \left |
\sum_{k=N}^{N_1}\sum_{s=M_0}^{m_0-1} c_{k,s}W_k(x) \cdot W_s(y)
\right | dx dy $$
$$+\max_{N<n \leq
N_1}\left[ \int\int_e \left | \sum_{k=N}^n  c_{k,m_0}W_k(x) \cdot
W_{m_0}(y) \right | dx dy\right] $$
$$\leq \left[ \int_{e_1}
\left | \sum_{k=N}^{N_1} a_kW_k(x) \right | dx \right] \cdot
\left[ \int_{e_2} \left | \sum_{s=M_0}^{m_0-1}  b_sW_s(y) \right |
dy \right]$$
$$+|b_{m_0}| \cdot  \left[ \int_{e_2} \left |
W_{m_0}(y) \right | dy \right] \cdot \max_{N<n \leq  N_1}\left[
\int_{e_1} \left | \sum_{k=N}^n a_kW_k(x) \right | dx \right] $$
$$\leq 12 \cdot |\gamma|\cdot |\Delta|.$$

Similarly, for $N \leq \overline {n} \leq N_1, \ \ M_0 \leq
\overline {m} \leq M$, we get
$$\int\int_e \left |
\sum_{k,s=N}^{\overline {n}, \overline{m}}  c_{k,s}W_k(x) \cdot
W_s(y) \right | dx dy  \leq 4\cdot |\gamma|\cdot |\Delta|.$$

{\bf Lemma 2 is proved.}
\par\par\bigskip

{\bf Lemma 3.} For any numbers $\varepsilon>0$, $N>1$ and a step
function
$$  {f(x,y)= \sum_{\nu=1}^{\nu_0} \gamma_\nu \cdot
\chi_{\Delta_\nu} (x,y)},$$ there exists a measurable set $E
\subset T$ and a polynomial $P(x,y)$ of the form
\[
P(x,y)= \sum_{k,s=N}^M c_{k,s}W_k(x) \cdot W_s(y) ,
\]
which satisfy the following conditions:
$$ P(x,y)=f(x,y) \ \  for \ \ (x,y) \in E,\leqno(I)$$
$$ |E|> 1- \varepsilon, \leqno(II)$$
$$\sum_{k,s=N}^M |c_{k,s}|^{2+\varepsilon}< \varepsilon,\leqno(III)$$
$$ \max_{N \leq \overline {n}, \overline{m} < M} \left[ \int\int_e
\left | \sum_{k,s=N}^{\overline {n}, \overline{m}}  c_{k,s}W_k(x)
\cdot W_s(y) \right | dx dy \right]\leqno(IV)$$
\[
+\max_{\sqrt {2}N \leq R \leq \sqrt {2}M} \left[ \int\int_e \left
| \sum_{2N^2 \leq k^2+s^2 \leq R^2}  c_{k,s}W_k(x) \cdot W_s(y)
\right | dx dy \right]
\]
\[
\leq 2 \cdot \int\int_e |f(x,y)| dxdy+ \varepsilon,
 \]
 for every measurable subset $e$  of $E$.
\vskip 4mm

{\bf Proof of Lemma 3.} Without any loss of generality, we
assume that
\begin{equation}
\max_{1 \leq \nu \leq \nu_0} \left( |\gamma_\nu| \cdot
|\Delta_\nu| \right)< {\varepsilon \over 32},
\end{equation}
($\Delta_ \nu,\ \ 1 \leq \nu \leq \nu_0$ are the constancy
rectangular domian of $f(x,y)$, i.e. where the function $f(x,y)$
is constant).

Given an integer $1 \leq \nu \leq \nu_0$, by applying Lemma 2 with $\delta={\varepsilon \over 16 \nu_0}$,
we find that there exists a measurable set $E_\nu \subset T$ and a
polynomial $P_\nu(x,y)$ of the form
\begin{equation}
P_\nu(x,y)= \sum_{k,s=N_\nu}^{M_\nu} c_{k,s}^{(\nu)}W_k(x) \cdot
W_s(y)
\end{equation}
with the following properties:
\begin{equation}
|E_\nu |> 1- {\varepsilon \over 2^\nu},
\end{equation}
\begin{equation}
\sum_{k,s=N_\nu}^{M_\nu} |c_{k,s}^{(\nu)}|^{2+\varepsilon}< {\varepsilon
\over \nu_0},
\end{equation}
\begin{equation}
P_\nu (x,y)=\gamma_\nu \cdot \chi_{\Delta_\nu} (x,y) \ \  for \ \
(x,y) \in E_\nu,
\end{equation}
$$\max_{N_\nu \leq \overline {n}, \overline{m} \leq M_\nu} \left[
\int\int_e \left | \sum_{k,s=N_\nu}^{\overline {n}, \overline{m}}
c_{k,s}^{(\nu)}W_k(x) \cdot W_s(y) \right | dx dy \right]$$
$$+\max_{\sqrt {2}N_\nu \leq R \leq \sqrt {2}M_\nu} \left[
\int\int_e \left | \sum_{2N_\nu^2 \leq k^2+s^2 \leq R^2}
c_{k,s}^{(\nu)}W_k(x) \cdot W_s(y) \right | dx dy \right] $$
\begin{equation}
\leq 16 \cdot |\gamma_\nu|\cdot |\Delta_\nu|<{\varepsilon \over 2},
\end{equation}
for every measurable subset $e$ of $E_\nu$ (see (3.5)).

Thus, for any fixed number $\nu,\ \ 1 \leq \nu \leq \nu_0$ we can
define a measurable set $E_\nu \subset T$ and a polynomial $P_\nu
(x,y)$ of the form (3.6), which satisfy the conditions (3.7) -
(3.9). Then we can take
$$N_1=N,\ \ n_\nu=M_{\nu-1}+1,\ \ 1 \leq \nu \leq \nu_0.$$
Set
\begin{equation}
E=\bigcap_{\nu=1}^{\nu_0} E_\nu,
\end{equation}
\begin{equation}
P(x,y)= \sum_{\nu=1}^{\nu_0}P_\nu(x,y)=\sum_{k,s=N}^M
c_{k,s}W_k(x) \cdot W_s(y),\ \ M=M_{\nu_0},
\end{equation}
where
\begin{equation}
c_{k,s}=c_{k,s}^{(\nu)},\ \ for \ \ N_\nu \leq k,s \leq M_\nu, \ \
1 \leq \nu \leq \nu_0
\end{equation}
and
\[
 c_{k,s}=0,\ \ for \ \ other \ \ k\ \ and\ \ s.
\]
>From (3.7) - (3.9), (3.11) - (3.13) we obtain:
$$P(x,y)=f(x,y) \ \  for \ \ (x,y) \in E, $$
$$|E|> 1- \varepsilon, $$
$$\sum_{k,s=N}^M |c_{k,s}|^{2+\varepsilon}< \sum_{\nu=1}^{\nu_0}
\left[ \sum_{k,s=N_\nu}^{M_\nu} |c_{k,s}^{(\nu)}|^{2+\varepsilon}
\right]< \varepsilon, $$ i.e the statements I) - III) of Lemma 3 are
satisfied. Now we will check the fulfillment of statement IV) of
Lemma 3.

Let $R \in [\sqrt {2}N,\sqrt {2}M ]$, then for some $\nu', 1 \leq
\nu' \leq \nu_0$ we have $\sqrt {2}N_{\nu'} \leq R \leq \sqrt
{2}N_{\nu'+1}$, consequently from (3.12) and (3.13) we have
\[
\sum_{2N^2 \leq k^2+s^2 \leq R^2}  c_{k,s}W_k(x) \cdot
W_s(y)=\sum_{\nu=1}^{\nu'-1}P_\nu(x,y)
 \]
 \[
 +\sum_{2N_{\nu'}^2 \leq
k^2+s^2 \leq R^2}  c_{k,s}^{(\nu')}W_k(x) \cdot W_s(y).
\]

In view of the conditions (3.9) - (3.12) and the equality $P(x,y)=f(x,y)$ on $E$,
for any measurable set $e \subset E$ we
obtain
\[
 \int\int_e \left | \sum_{2N^2 \leq k^2+s^2 \leq R^2}
c_{k,s}W_k(x) \cdot W_s(y) \right | dxdy
\]
\[
\leq \int\int_e \left | \sum_{\nu=1}^{\nu'-1}P_\nu(x,y) \right |
dxdy
\]
\[
+ \int\int_e \left | \sum_{2N_{\nu'}^2 \leq k^2+s^2 \leq R^2}
c_{k,s}^{(\nu')}W_k(x) \cdot W_s(y) \right | dxdy
\]
\[
\leq \int\int_e |f(x,y)|dxdy+{\varepsilon \over 2}.
\]
Similarly, for any $e \subset E$ we have
\[
\max_{N \leq \overline {n}, \overline{m} \leq M} \left[ \int\int_e
\left | \sum_{k,s=N}^{\overline {n}, \overline{m}}  c_{k,s}W_k(x)
\cdot W_s(y) \right | dx dy \right]
\]
\[
 \leq \int\int_e
|f(x,y)|dxdy+{\varepsilon \over 2}.
\]

{\bf Lemma 3 is proved.}
\par\par\bigskip

\section{Proofs of the theorems}

\vskip 4mm

{\bf Proof of Theorem 3.}
 \vskip 4mm Let
\begin{equation}
\{ f_s(x,y)\}_{s=1}^\infty,\ \  (x,y) \in T
\end{equation}
 be a sequence of all step functions, values and constancy interval endpoints of which
are rational numbers. Applying Lemma 3 consecutively, we can find
a sequence
  $\{ E_s\}_{s=1}^\infty $ of sets and a sequence of polynomials
\begin{equation}
P_s(x,y)=\sum_{k,\nu=N_{s-1}}^{
N_s-1}c_{k,\nu}^{(s)}W_k(x)W_\nu(y) ,
\end{equation}
\[
 1=N_0<N_1<...<N_s<....,\ \ s=1,2,....,
\]
which satisfy the conditions:
\begin{equation}
P_s(x,y)=f_s(x,y),\ \ (x,y)\in E_s,
\end{equation}
\begin{equation}
\left| E_s\right| >1-2^{-2(s+1)} ,\ \  E_s\subset T,
\end{equation}
\begin{equation}
\sum_{k,\nu=N_{s-1}}^{ N_s-1}\left
|c_{k,\nu}^{(s)}\right|^{2+2^{-2s}}< 2^{-2s},
\end{equation}
\[
\max_{N_{s-1} \leq \overline {n}, \overline{m} < N_{s}} \left[
\int\int_e \left | \sum_{k,\nu=N_{s-1}}^{\overline {n},
\overline{m}} c_{k,\nu}^{(s)}W_k(x) \cdot W_\nu(y) \right | dx dy
\right]
\]
\[+\max_{\sqrt {2}N_{s-1} \leq R \leq \sqrt {2}N_s} \left[ \int\int_e \left
| \sum_{2N_{s-1}^2 \leq k^2+\nu^2 \leq R^2}  c_{k,\nu}^{(s)}W_k(x)
\cdot W_\nu(y) \right | dx dy \right]
\]
\begin{equation}
\leq 2 \cdot \int\int_e |f_s(x,y)| dxdy+ 2^{-2(s+1)},
\end{equation}
for every measurable subset $e$  of $E_s$.

Denote
\begin{equation}
\sum_{k,\nu=1}^\infty c_{k,\nu} W_k(x)W_\nu(y)=\sum_{s=1}^\infty
\left[\sum_{k,\nu=N_{s-1}}^{ N_s-1}c_{k,\nu}^{(s)}W_k(x)W_\nu(y)
\right],
\end{equation}
where
\[
c_{k,\nu}=c_{k,\nu}^{(s)}, \ \  for\ \ N_{s-1}\leq k,\nu <N_s,\ \
s=1,2,... .
\]

For an arbitrary number $\varepsilon >0$  we set
\[
 \Omega_n =\displaystyle{ \bigcap_{s=n}^\infty E_s,}\ \  n=1,2,....;\hfil
\quad
\]
\begin{equation}
 E=\Omega_{n_0} =\displaystyle{ \bigcap_{s=n_0}^\infty E_s,}\ \   n_0=[\log_{1/2} \varepsilon]+1; \hfil
 \end{equation}
 \[
 B=\displaystyle{\bigcup _{n=n_0} ^\infty \Omega_n} =\Omega_{n_0}
\bigcup \left( \displaystyle{\bigcup _{n=n_0+1}^ \infty }\Omega_n
\setminus \Omega_{n-1} \right)
\]

It is obvious ( see (4.4), (4.8) ) that $\left| B \right|=1$  and
$\left| E \right| >1- \varepsilon .$

We define a function $\mu(x,y)$ in the following way:
\begin{equation}
\mu(x,y)=1, \ \ for \ \ (x,y) \in E \cup (T \setminus B);
\end{equation}
\[
 \mu(x,y)= \mu_n, \ \ for  \ \ (x,y) \in \Omega_n \setminus \Omega_{n-1},\ \ n\geq
n_0+1,
\]
where
\begin{equation}
 \displaystyle{ \mu_n=\left[
2^{2n}\cdot \prod_{s=1}^n h_s \right]^{-1}};\hfil \quad
\end{equation}
\[
 h_s=|| f_s||_C+\displaystyle{ \max_{N_{s-1} \leq \overline {n}, \overline{m} < N_{s}}
 \left|\left|  \sum_{k,\nu=N_{s-1}}^{\overline {n},
\overline{m}} c_{k,\nu}^{(s)}W_k(x) \cdot W_\nu(y)\right|
\right|_C} \hfil
\]
\[
+ \displaystyle{ \max_{\sqrt {2}N_{s-1} \leq R \leq \sqrt {2}N_s}
\left|\left| \sum_{2N_{s-1}^2 \leq k^2+\nu^2 \leq R^2}
c_{k,\nu}^{(s)}W_k(x) \cdot W_\nu(y)\right| \right|_C +1.}
\]

>From (4.5), (4.7) - (4.10) we obtain

$(A) - 0<\mu(x,y) \le1, \ \ \mu(x,y)$ is a measurable function and
\[
\left | \{(x,y)\in T:\mu(x,y)\not =1\} \right|<\varepsilon.
\]
(B) -- $\displaystyle {\sum_{k,\nu=1}^\infty
\left|c_{k,\nu}\right|^q<\infty} \ \mbox{for all}\ q>2.$

Hence, obviously we have (see (4.5) and (4.7))
$$ \lim_{k,\nu\to \infty}c_{k,\nu}=0. $$
It follows  from (4.8) - (4.10) that for all $s \geq n_0$  and
$N_{s-1} \leq \overline {n}, \overline{m} < N_{s}$
\[
\int\int_{T \setminus \Omega_s} \left|
\sum_{k,\nu=N_{s-1}}^{\overline {n}, \overline{m}}
c_{k,\nu}^{(s)}W_k(x) \cdot W_\nu(y)\right| \mu(x,y) dxdy
\]
\[
=\sum_{n=s+1}^ \infty \left[\int\int_{\Omega_n \setminus
\Omega_{n-1}} \left| \sum_{k,\nu=N_{s-1}}^{\overline {n},
\overline{m}} c_{k,\nu}^{(s)}W_k(x) \cdot W_\nu(y)\right| \mu_n
dxdy \right]
\]
\begin{equation}
\leq \sum_{n=s+1}^ \infty2^{-2n} \left[\int\int_T \left|
\sum_{k,\nu=N_{s-1}}^{\overline {n}, \overline{m}}
c_{k,\nu}^{(s)}W_k(x) \cdot W_\nu(y)\right| h_s^{-1} dxdy
\right]<{1 \over 3}2^{-2s}.
\end{equation}
Analogously for all $s \geq n_0$  and $\sqrt {2}N_{s-1} \leq R
\leq \sqrt {2}N_s$ we have
\begin{equation}
\int\int_{T \setminus \Omega_s} \left|\sum_{2N_{s-1}^2 \leq
k^2+\nu^2 \leq R^2} c_{k,\nu}^{(s)}W_k(x) \cdot W_\nu(y)\right|
\mu(x,y) dxdy<{1 \over 3}2^{-2s}.
\end{equation}
By (4.2), (4.8) - (4.10) for all  $s \geq n_0$ we have
\[
\int\int_T \left| P_s(x,y)-f_s(x,y) \right|\mu(x,y)dxdy
\]
\[
=\int\int_{\Omega_s} \left| P_s(x,y)-f_s(x,y) \right|\mu(x,y)dxdy
\]
\[+\int\int_{T \setminus \Omega_{s}} \left| P_s(x,y)-f_s(x,y)
\right|\mu(x,y)dxdy
\]
\[
=\sum_{n=s+1}^\infty \left[\int\int_{\Omega_n \setminus
\Omega_{n-1}} \left| P_s(x,y)-f_s(x,y)  \right| \mu_n dxdy\right]
\]
\[
\leq \sum_{n=s+1}^ \infty 2^{-2n}\left[ \int\int_T \left(\left|
f_s(x,y) \right| + \sum_{k,\nu=N_{s-1}}^ {N_s-1}
c_{k,\nu}^{(s)}W_k(x) \cdot W_\nu(y) \right) h_s^{-1}dxdy \right ]
\]
\begin{equation}
<{1 \over 3}2^{-2s}<2^{-2s}.
\end{equation}
By (4.6) and (4.8) - (4.11) for all
 $N_{s-1}\leq \overline {n}, \overline{m} < N_{s}$ and $s
\geq n_0+1$ we obtain
\[
\int\int_T  \left| \sum_{k,\nu=N_{s-1}}^{\overline {n},
\overline{m}} c_{k,\nu}^{(s)}W_k(x) \cdot W_\nu(y)\right| \mu(x,y)
dxdy
\]
\[
\int\int_{\Omega_s}  \left| \sum_{k,\nu=N_{s-1}}^{\overline {n},
\overline{m}} c_{k,\nu}^{(s)}W_k(x) \cdot W_\nu(y)\right| \mu(x,y)
dxdy
\]
\[
\int\int_{T \setminus \Omega_s}  \left|
\sum_{k,\nu=N_{s-1}}^{\overline {n}, \overline{m}}
c_{k,\nu}^{(s)}W_k(x) \cdot W_\nu(y)\right| \mu(x,y) dxdy
\]
\[
< \sum_{n=n_0+1}^ s \left[\int\int_{\Omega_n \setminus
\Omega_{n-1}} \left| \sum_{k,\nu=N_{s-1}}^{\overline {n},
\overline{m}} c_{k,\nu}^{(s)}W_k(x) \cdot W_\nu(y)\right|\cdot
\mu_ndxdy \right]+{1 \over 3}2^{-2s}
\]
\[
<\sum_{n=n_0+1}^ s \left( 2^{-2(s+1)}+2\cdot\int\int_{\Omega_n
\setminus \Omega_{n-1}} |f_s(x,y)|dxdy \right)\cdot \mu_n +{1
\over 3}2^{-2s}
 \]

\[
=2^{-2(s+1)} \cdot \sum_{n=n_0+1}^ s \mu_n+
+\int \int_{\Omega_s} |f_s(x,y)|\mu(x,y)dxdy +{1 \over 3}2^{-2s}
\]
\begin{equation}
<2\cdot\int \int_T |f_s(x,y)|\mu(x,y)dxdy +2^{-2s}.
\end{equation}
Analogously for all $s \geq n_0$  and $\sqrt {2}N_{s-1} \leq R
\leq \sqrt {2}N_s$ we have (see (4.12))
\[
\int\int_T \left|\sum_{2N_{s-1}^2 \leq k^2+\nu^2 \leq R^2}
c_{k,\nu}^{(s)}W_k(x) \cdot W_\nu(y)\right| \mu(x,y) dxdy
\]
\begin{equation}
<2\cdot\int \int_T |f_s(x,y)|\mu(x,y)dxdy +2^{-2s}.
\end{equation}
Now we'll show that the series (4.7) is universal in $L_\mu^1(T)$
concerning subseries with respect to convergence by both spherical
and rectangular partial sums.

Let $ f(x,y) \in L_{\mu}^1 (T)$ , i. e.
\[
\int\int_T |f(x,y)|\mu(x,y) dxdy<\infty.
\]

It is easy to see that we can choose a function $f_{n_1}(x,y)$
from the sequence (4.1) such that
\begin{equation}
\int\int_T \left| f(x,y)- f_{n_1}(x,y)
\right|\mu(x,y)dxdy<2^{-2},\ \ n_1
> n_0+1.
\end{equation}
 Hence, we have
\begin{equation}
\int\int_T \left| f_{n_1}(x,y)
\right|\mu(x,y)dxdy<2^{-2}+\int\int_T |f(x,y)|\mu(x,y)dxdy.
 \end{equation}
 From (4.13) and (4.16) we get
\[
\int\int_T \left| f(x,y)- P_{n_1}(x,y) \right|\mu(x,y)dxdy
\]
\[
\leq \int\int_T \left| f(x,y)- f_{n_1}(x,y) \right|\mu(x,y)dxdy
\]
\[
+\int\int_T \left| f_{n_1}(x,y)- P_{n_1}(x,y) \right|\mu(x,y)dxdy<
2\cdot 2^{-2}.
 \]
Assume that numbers $n_1<n_2<...<n_{q-1}$ are chosen in such a way
that the following condition is satisfied:
\begin{equation}
\int\int_T \left| f(x,y)- \sum_{s=1}^j P_{n_s}(x,y)
\right|\mu(x,y)dxdy<2\cdot 2^{-2j}, \ \ 1\leq j \leq q-1 .
\end{equation}
Now we choose a function $f_{n_q}(x,y)$ from the sequence (4.1)
such that
\[
\int\int_T \left| \left( f(x,y)- \sum_{s=1}^{q-1} P_{n_s}(x,y)
\right)-f_{n_q}(x,y)\right| \mu(x,y)dxdy
\]
\begin{equation}
 <2\cdot 2^{-2q},\ \ n_q>n_{q-1}.
\end{equation}
This with (4.18) imply
\begin{equation}
\int\int_T \left| f_{n_q}(x,y)\right| \mu(x,y)dxdy<
 2^{-2q}+2\cdot 2^{-2(q-1)}=9\cdot 2^{-2q}.
\end{equation}
 From (4.2), (4.13) - (4.15) and (4.20) we obtain
\begin{equation}
\int\int_T \left| f_{n_q}(x,y)- P_{n_q}(x,y) \right|\mu(x,y)dxdy<
 2^{-2n_q},
\end{equation}
where
\[
P_{n_q}(x,y)=\sum_{k,\nu=N_{n_q-1}}^{
N_{n_q}-1}c_{k,\nu}^{(n_q)}W_k(x)W_\nu(y) ,
\]
\[
\max_{N_{n_q-1} \leq \overline {n}, \overline{m}<N_{n_q}} \left[
\int\int_T \left| \sum_{k,\nu=N_{n_q-1}}^{\overline {n},
\overline{m}} c_{k,\nu}^{(n_q)}W_k(x) \cdot W_\nu(y)\right|
\mu(x,y) dxdy \right ]
\]
\begin{equation}
<19 \cdot 2^{-2q}.
\end{equation}
Analogously we have
\[
\max_{\sqrt {2}N_{n_q-1} \leq R \leq \sqrt {2}N_{n_q}} \left[
\int\int_T \left|\sum_{2N_{n_q-1}^2 \leq k^2+\nu^2 \leq R^2}
c_{k,\nu}^{(n_q)}W_k(x) \cdot W_\nu(y)\right| \mu(x,y) dxdy\right]
\]
\[
<19 \cdot 2^{-2q}.
\]
 In quality subseries of Theorem we shall take
\begin{equation}
\sum_{q=1}^\infty P_{n_q}(x,y)=\sum_{q=1}^\infty
\left[\sum_{k,\nu=N_{n_q-1}}^{
N_{n_q}-1}c_{k,\nu}^{(n_q)}W_k(x)W_\nu(y) \right].
\end{equation}
>From (4.19) and (4.21) we have
\[
\int\int_T \left| f(x,y)-\sum_{s=1}^q P_{n_s}(x,y)\right| \mu(x,y)
dxdy
\]
\[
\leq \int\int_T \left|\left( f(x,y)- \sum_{s=1}^{q-1} P_{n_s}(x,y)
\right)-f_{n_q}(x,y)\right| \mu(x,y)dxdy
\]
\begin{equation}
+\int\int_T \left| f_{n_q}(x,y)- P_{n_q}(x,y) \right|\mu(x,y)dxdy<
 2\cdot 2^{-2q}.
\end{equation}

Let $ \overline {n}$ and $\overline{m}$ be arbitrary natural
numbers. Then for some natural number $q$ we have
\[
 N_{n_q-1}\leq \min \{\overline {n}, \overline{m}\}<N_{n_q}.
\]
Taking into account (4.22) and (4.24) for rectangular partial sums
$S_{\overline {n},\overline {m}}(x,y)$ of (4.23) we get
\[
\int\int_T \left| S_{\overline {n},\overline
{m}}(x,y)-f(x,y)\right| \mu(x,y) dxdy
\]
\[\leq \int\int_T \left|
f(x,y)-\sum_{s=1}^q P_{n_s}(x,y)\right| \mu(x,y) dxdy
\]
\[
+\max_{N_{n_q-1} \leq \overline {n}, \overline{m}<N_{n_q}} \left[
\int\int_T \left| \sum_{k,\nu=N_{n_q-1}}^{\overline {n},
\overline{m}} c_{k,\nu}^{(n_q)}W_k(x) \cdot W_\nu(y)\right|
\mu(x,y) dxdy \right ]
\]
\begin{equation}
 <21 \cdot 2^{-2q}.
\end{equation}
Analogously for $\sqrt {2}N_{n_q-1} \leq R \leq \sqrt {2}N_{n_q}$
we have
\begin{equation}
\int\int_T \left|S_R(x,y)-f(x,y)\right| \mu(x,y) dxdy<21 \cdot
2^{-2q},
\end{equation}
where $S_R(x,y)$ the spherical partial sums of (4.23).

>From (4.25) and (4.26) we conclude that the series (4.7) is
universal in $L_\mu^1(T)$ concerning subseries with respect to
convergence by both spherical and rectangular partial sums (see
Definition 2).

\vskip 4mm {\bf Theorem 3 is proved.}

\vskip 4mm

 {\bf Remark.} We can show Theorem 4 by the method in the proof of Theorem 3.

\vskip 4mm \break

\vspace{5mm}

\author{Department of Physics, State University of Yerevan,
       Alex Manukian 1, 375049 Yerevan, Armenia}

{\em E-mail address: sergoep@ysu.am}


\begin{thebibliography}{10}

\vskip 4mm

\bibitem{Arut76} F.G.Arutunian, To representation of functions by multiplaer series, Dokl. Akad. Nauk Armyan. SSR 64(1976), 72-76 [in russian].
\bibitem{Dza64} O.P.Dzangadze, On the universal double series, Bull. Georgian Acad. Sci. 34(1964), 225-228 [in russian].
\bibitem{Epis99} S.A. Episkoposian , "On the series by Walsh system universal in weighted $L_\mu^1[0,1]$ spaces ",
  Izv. Nats. Akad. Nauk Armenii, Math., English trans. in: J. Contemp. Math. Anal., 34(1999), no. 2, 25-40.
\bibitem{Epis06} S.A. Episkoposian , "On the existence of universal series by trigonometric system ",
J. Funct. Anal., 230(2006), 169 - 189.
\bibitem{Grig99} M.G.Grigorian "On the representation of functions by orthogonal series in weighted $L^p$ spaces,
\ \ Studia Math. 134(1999), 211-237.
\bibitem{Grig01} M.G.Grigorian, S.A. Episkoposian  "Representation of functions
in weighted spaces $L_\mu^1[0,1]$ by  trigonometric and Walsh
series, \ \ Anal. Math., 27 (2001), 261-277.
\bibitem{Kozl50} V.Ya.Kozlov, On the complete systems of orthogonal functions, Mat. Sb.  26(1950), 351-364 [in russian].
\bibitem{Mensh47} D.E.Menshov, On the partial summs of trigonometric series, Mat. Sb. 20(1947), 197-238 [in russian].
\bibitem{Orl27} W.Orlicz , Uber die unabhangig von der Anordnung fast
     uberall kniwergenten Reihen, Bull. Acad. Polonaise Sci., 81 (1927), 117-125.
\bibitem{Pal32} R.Paley, A remarkable systems of orthogonal functions, Proc. London
Math. Soc. 34(1932), 241-279.
\bibitem{Tala57} A.A.Talalian, On the convergence almost everywhere the subsequence
of partial sums of general orthogonal series, Izv. Akad. Nauk
Armyan. SSR Ser. Math. 10(1957), 17-34 [in russian].
\bibitem{Tala60} A.A.Talalian, On the universal series with respect to rearrangements,
Izv. Akad. Nauk SSSR Ser. Math. 24(1960), 567-604 [in russian].


\end{thebibliography}
\end{document}